\documentclass[12pt,fleqn]{article} 
\usepackage{amsfonts,amssymb}
\newcommand{\ds}{\displaystyle}
\newcommand{\y}{\\[3pt]}
\newcommand{\q}{\quad}\newcommand{\qq}{\qquad}
\newcommand{\tw}{\tilde\w} 
\renewcommand\ge{\geqslant}\renewcommand\le{\leqslant}
\newcommand\bb{\bigbreak}\newcommand\n{\noindent}
\newcommand\mb{\medbreak}\renewcommand\sb{\smallbreak}
\newcommand{\g}{\mathfrak{g}}
\newcommand{\cF}{\mathcal{F}}

\newcommand{\ad}{\mathrm{ad}}\newcommand{\tr}{\mathrm{tr}}
\newcommand{\C}{\mathbb{C}}\newcommand{\R}{\mathbb{R}}
\newcommand{\CP}{\mathbb{CP}}

\newcommand{\alr}{\al\,.\,,.\,\ar}\renewcommand{\=}{&\!=\!&\ds}
\newcommand{\Ric}{\mathrm{Ric}}
\newcommand{\cir}{\hbox{\small$\circ$}}
\renewcommand{\a}{\alpha}\renewcommand{\b}{\beta}
\newcommand{\fa}{\mathfrak{a}}\newcommand{\fb}{\mathfrak{b}}
\newcommand{\fh}{\mathfrak{h}}
\newcommand{\fr}{\mathfrak{r}}\newcommand{\fs}{\mathfrak{s}}
\newcommand{\fu}{\mathfrak{u}}\newcommand{\fz}{\mathfrak{z}}
\newcommand{\cR}{\mathcal{R}}\newcommand{\tx}{{\tilde x}}
\newcommand{\cV}{\mathcal{V}}\newcommand{\cW}{\mathcal{W}}
\newcommand{\der}{\hbox{der}\,}\renewcommand{\Im}{\hbox{Im}}
\renewcommand{\t}{\lambda}\newcommand{\s}{k}
\renewcommand{\th}{\theta}\newcommand\w{\omega}
\newcommand{\e}{\mathbf{e}}\newcommand{\f}{\mathbf{f}}

\newcommand{\ph}{\vphantom{\int^l}}\newcommand{\we}{\wedge}

\newcommand{\op}{\oplus}\newcommand{\ot}{\otimes}
\newcommand\qed{\hfill$\square$\mb}\newcommand{\rf}[1]{(\ref{#1})}
\newcommand{\E}{\raise1pt\hbox{\small$\textstyle\bigwedge$}\kern-0.5pt}
\newcommand{\ba}{\begin{array}}\newcommand{\ea}{\end{array}}
\newcommand{\be}{\begin{equation}}
\newcommand{\ee}[1]{\label{#1}\end{equation}}
\newcommand{\bt}{\begin{tabular}}\newcommand{\et}{\end{tabular}}
\newcommand{\al}{\langle}\newcommand{\ar}{\rangle}
\newcommand{\ft}[2]{\hbox{$\textstyle\frac{#1}{#2}$}}
\newcommand{\frs}[2]{\hbox{\large$\textstyle\frac{#1}{#2}$\normalsize}}
\newcommand{\nbf}[1]{\stepcounter{equation}
                                         \bb\n\textbf{\theenumii\ #1.}}
\newcommand{\nit}[1]{\sb\n\textit{#1.}}
\renewcommand{\theenumii}{\arabic{subsection}.\arabic{equation}}

\begin{document}\parskip3pt
\title{\bf Anti-self-dual metrics on Lie groups}\vspace{20pt}
\author{Vivian De Smedt and Simon Salamon}
\date{}
\maketitle
 
\nit{Abstract} The aim of the paper is to determine left-invariant,
anti-self-dual, non conformally flat, Riemannian metrics on four-dimensional
Lie groups.\mb

\setcounter{equation}0
\subsection*{Introduction}

Let $G$ be a Lie group of dimension $n$ equipped with a left-invariant
Riemannian metric $g$. The curvature tensor $R=R(g)$ is completely determined
by its value at the identity element of $G$. In this way, $R$ is an element of
the vector space \[\cR=S^2(\E^2\g^*)\ominus\E^4\g^*,\] where $\g$ is the Lie
algebra of $G$ and $A\ominus B$ denotes the orthogonal complement of $B$ in
$A$. Moreover, there is an $O(n)$-invariant decomposition
\[\cR=\cV\op\cW,\] where $\cV\cong S^2\g^*$. The component $\Ric$ of $R$ in 
$\cV$ is the Ricci tensor, and the component $W$ in $\cW$ is the Weyl
tensor. If $W=0$, then $M$ is conformally flat in the sense that there exist
local coordinates $x^i$ for which $g$ is a scalar function times $\sum_{i=1}^n
dx^i\otimes dx^i$. For example, the compact Lie group $U(2)$ (double covered by
$S^1\times S^3$) acquires such coordinates from its description as a discrete
quotient of $\R^4\setminus\{0\}$.

If $G$ is 4-dimensional and oriented, the Hodge involution $*\colon\E^2\g^*\to
\E^2\g^*$ gives rise to a decompostion \[W=W^++W^-,\] where $W^\pm\in
S^2_0(\E^\pm\g^*)$ and $\E^\pm\g^*$ is the $\pm$-eigenspace of $*$. If $W^+=0$
the metric $g$ is said to be \textit{anti-self-dual}, or `ASD' for short. In
this paper, we classify left-invariant metrics on 4-dimensional Lie groups with
$W^+=0$ and $W^-\ne0$. Such a group is necessarily solvable, and in \S3 we
describe the general form of its Lie algebra. We prove that the hypotheses
eliminate all but a 1-parameter family of Lie algebras, and that there are
essentially only two distinct left-invariant metrics that are anti-self-dual.

The proof of the main theorem ({\bf1.6} below) is accomplished by first
eliminating Lie algebras that admit an orientation-reversing automorphism. A
key feature of the method is that the Gram-Schmidt process is compatible with
the triangular nature of the nilpotent Lie algebra structure constants. This
enables one to compute the structure constants for an orthonormal basis
relative to an arbitrary inner product. Once freedom in the choice of basis
has been removed, the equation $W^+=0$ is solved by means of a specially
adapted Maple program reproduced in the Appendix. No doubt further effort would
eliminate the need for a computer analysis, but this aspect of the work may
have independent interest.\sb

\n{Acknowledgments.} This paper was begun during a visit by the first author
in Oxford in 1997, in parallel with the work \cite{Abb} that establishes
similar results independently. The second author is grateful for the
encouragement received at the Tokyo conference to finally present the material
in written form.

\setcounter{equation}0\setcounter{subsection}1
\subsection*{1. Summary of results}

Let $V$ denote a real 4-dimensional vector space, with basis $(\f_1,\f_2,\f_3,
\f_4)$. For $\t\in\R$, let $\g_\t$ denote the Lie algebra defined on $V$ by the
relations \be\ba{l}[\f_1,\f_2]=\f_2-\t \f_3,\y[\f_1,\f_3]=\t\f_2+\f_3,\y
[\f_1,\f_4]=\>2\f_4, \y[\f_2,\f_3]=-\f_4.\ea\ee{LAt} The Jacobi identity is
easily verified. The Lie algebra $\g_\t$ is solvable and $\g_\t'=[\g_\t,\g_\t]$
is isomorphic to the 3-dimensional Heisenberg algebra. Changing the signs of
$\f_3,\f_4$ is equivalent to replacing $\t$ by $-\t$, though we shall see that
$\g_\t,\,\g_{\t'}$ are not isomorphic if $|\t|\ne|\t'|$.

Later on, it will be convenient to express the relations \rf{LAt} in terms of
differential forms. To this end, let $(\f^i)$ denote the dual basis of
$\g_\t^*$. With the convention that $d\f^i(u,v)=-\f^i[u,v]$, \rf{LAt} asserts
that \be\ba{l} d\f^1=0,\y d\f^2=-\f^{12}-\t \f^{13},\y d\f^3=\>\t
\f^{12}-\f^{13},\y d\f^4=-2\f^{14}+\f^{23}.\ea\ee{dfi} Let $G_\t$ denote the
simply-connected Lie group corresponding to $\g_\t$. We shall not identify the
group structure, but \rf{dfi} can be integrated to yield coordinates $t,x,y,u$
on $G_\t$ relative to which \be\ba{l} \f^1=dt,\y \f^2=e^{-t}(dx+\t ydt),\y
\f^3=e^{-t}(dy-\t xdt),\y
\f^4=e^{-2t}(du+\ft12(xdy-ydx)-\ft12\t(x^2\!+\!y^2)dt).\ea\ee{1f}

For $\s>0$, let $\alr_\s$ denote the inner product on $V$ defined by \be\al
\f_i,\f_j\ar=\left\{\ba{ll}\s^2,\qq& i=j=1,\y 1,& 2\le i=j\le4,\y 0,& i\ne
j.\ea\right.\ee{ip} This inner product induces a left-invariant metric
\[ g_\s=\s^2 \f^1\ot \f^1 +\f^2\ot \f^2+\f^3\ot \f^3+\f^4\ot \f^4\] on $G_\t$. 
We shall show that the corresponding metrics are all isometric for fixed
$\s$. The associated curvature tensor does not therefore depend upon
$\t$. Using the program in the Appendix, it is easy to verify

\nbf{Proposition} Relative to the orthonormal basis $(\frs1\s
\f_1,\f_2,\f_3,\f_4)$ of $\g_\t$, the Ricci and Weyl tensors of $g_\s$ have the
diagonal forms

\[\Ric=\left(\ba{cccc}\frs6{\s^2}&0&0&0\\0&\frs4{\s^2}+\frac12&0&0\\
0&0&\frs4{\s^2}+\frac12&0\\0&0&0&\frs8{\s^2}-\frac12\ea\right),\]
\[W^\pm=\frac{\s^2\mp 3\s +2}{3\s^2}\!\left(\ba{ccc}-1&0&0\\0&-1&0\\
0&0&2\ea\right).\] 

Observe that $\Ric$ is degenerate if $\s=4$, and $W^+=0$ if $\s\in\{1,2\}$. We
shall prove that, conformally speaking, there are no other non-flat ASD
examples:

\nbf{Theorem} Let $G$ be an oriented four-dimensional Lie group admitting a
left-invariant Riemannian metric $g$ such that $W^+=0$ and $W^-\ne0$. Then the
Lie algebra of $G$ is isomorphic to $\g_\t$ for some $\t\ge0$, and $g$ is
locally homothetic to either $g_1$ or $g_2$.\mb

The distinction between $W^+,W^-$ depends upon a choice of orientation, and in
this sense inclusion of the prefix `anti' in the title is purely a matter of
taste. However, an orientation \textit{is} distinguished by the presence of a
complex structure. The set of positively-oriented almost complex structures on
$G$ compatible with a metric is isomorphic to $SO(4)/U(2)\cong S^2$, and the
equation $W^+=0$ is precisely the integrability condition for a tautological
almost complex structure on $G\times S^2$ \cite{AHS}. A special situation in
which $W^+=0$ occurs when each element in $S^2$ represents an integrable
complex structure on $G$, which is then called \textit{hypercomplex}.

Four-dimensional Lie groups with a left-invariant hypercomplex structure were
classified by Barberis \cite{Bar}, and the above theorem can therefore be
viewed as a generalization of this work. Indeed, \cite{Bar} asserts that
$G_\t$ admits a left-invariant hypercomplex structure if and only if $\t=0$,
and that $(G_0,g_1)$ is hyperhermitian. On the other hand, $(G_0,g_2)$ is
isometric to the complex hyperbolic plane $\mathbb{CH}^2=SU(2,1)/S(U(2)\times
U(1))$ with its symmetric metric, which is well known to satisfy $W^-=0$ (with
respect to the orientation relative to the natural complex structure). This
metric appeared in Jensen's classification of Einstein metrics on Lie groups
\cite{Jen}, which provides a starting point for an alternative approach to
classifying self-dual metrics and other structures on 4-dimensional Lie groups
\cite{Abb}.

The fact that there are just two values of the parameter $\s$ that solve our
problem is reminiscent of the existence of two Einstein metrics on certain
sphere bundles. For example, it is well known that $\CP^3$ (as a bundle over
$S^4$) has its standard K\"ahler-Einstein metric as well as a nearly-K\"ahler
metric with weak holonomy $U(3)$.

\setcounter{equation}0\setcounter{subsection}2
\subsection*{2. Building Lie algebras}

Let $\fa_k$ denote the abelian Lie algebra whose underlying vector space is
$\R^k$. Let $\fh_3$ denote the Lie algebra of the Heisenberg group; thus
$\fh_3$ has a basis $(\f_1,\f_2,\f_3)$ satisfying $\f_1=[\f_2,\f_3]=-[\f_3,\f_2]$ and
all other brackets zero.

Given a Lie algebra $\fh$, let \[\der\fh=\{D\colon\fh\to\fh:D[x,y]=[Dx,y]+
[x,Dy]\}\] denote the set of derivations of $\fh$. Suppose that
$\rho\colon\fb\to\der\fh$ is a homomorphism of Lie algebras. The
\textit{extension of $\fh$ by $\rho$} is the Lie algebra $\g$ with underlying
vector space $\fb\op\fh$, in which $\fb$ and $\fh$ are both subalgebras and
$[x,h]=\rho(x)(h)$ for all $x\in\fb$ and $h\in\fh$. We write
\[ \g=\fb\op_\rho\fh,\] and if $\fb$ is abelian, we say that $\g$ is an 
\textit{abelian extension} of $\fh$.

\nbf{Proposition} If $\g$ is a 4-dimensional Lie algebra with zero centre then
$\g$ is isomorphic to an abelian extension of $\fa_2$, $\fa_3$ or $\fh_3$.

\nit{Proof} If $\g$ is not solvable, its radical $\fr$ is a proper ideal. The
Levi decomposition gives $\g=\fs\op_\rho\fr$ where $\fs$ is a semi-simple
algebra and $\rho\colon\fs\to\der\fr$. It follows that $\dim\fs=3$ and
$\dim\fr=1$. But 1-dimensional representations of a semi-simple Lie algebra are
trivial, so the centre $\fz$ of $\g$ equals $\fr$. We may therefore assume that
$\g$ is solvable.

A Lie algebra $\g$ is solvable if and only if the ideal $\g'=[\g,\g]$ is
nilpotent. The only non-abelian nilpotent Lie algebra of dimension less than 4
is $\fh_3$, so $\g'$ is one of $0$, $\fa_1$, $\fa_2$, $\fa_3$, $\fh_3$.\sb

\n If $\g'=0$ then $\fz=\g\ne0$.\sb 

\n If $\g'\cong\fa_1$ then $\g$ has a basis $(\f_1,\f_2,\f_3,\f_4)$ with $\g'=\al
\f_1\ar$. Unless $\f_1\in\fz$, we may suppose that $\f_2,\f_3,\f_4$ are chosen such
that $[\f_1,\f_2]=\f_1$ and $[\f_1,\f_i]=0$ for $i=3,4$.  We may further modify
$\f_3$ in order that $[\f_2,\f_3]=0$. But then the Jacobi identity \[[[\f_2,\f_3],
\f_4]+[[\f_3,\f_4],\f_2]+[[\f_4,\f_2],\f_3]=0\] implies that $[\f_3,\f_4]=0$ and
$\f_3\in\fz$.\sb

\n If $\g'\cong\fa_2$ then there is an isomorphism $\g\cong\R^2\op\g'$ of
vector spaces, and the bracket determines a linear mapping $\rho\colon\R^2\to
\der(\g')$. If $\ker\rho$ is non-zero then it contains a 1-dimensional subspace
that when added to $\g'$ yields an abelian algebra $\fa_3$ of which $\g$ is an
abelian extension. If $\rho$ is injective then $\g'_\C$ has a basis of
eigenvectors for the commuting elements in $\Im\rho$. It follows that there
exists a real basis $(\f_1,\f_2,\f_3,\f_4)$ of $\g$ such that $\g'=\al \f_3,\f_4\ar$
and $\ad \f_1$ restricts to the identity on $\g'$. Then $\f_2$ may be modified so
that $[\f_1,\f_2]=0$, and $\g=\R^2\op_\rho\g'$ is an abelian extension.

\n If $\g'\cong\fa_3$ or $\g'\cong\fh_3$ then $\g$ is immediately an abelian
extension of $\g'$.\qed

Four-dimensional solvable Lie algebras can be broadly divided into 7 classes
according to the triple $(d',d'',d''')$ of dimensions of $\g'$,
$\g''\!=\![\g',\g']$, $\g'''\!=\![\g'',\g'']$.  In \cite{Abb}, curvature
computations are carried out for each of these classes in turn. For our
purposes, the following observation helps to restrict the range of algebras
that need to be considered.

\nbf{Lemma} Let $G$ be a simply-connected Lie group with a left-invariant
metric $g$. There exists an orientation-reversing isometry of $(G,g)$ in either
of the two cases:\\\n(i) the Lie algebra $\g$ has non-zero centre, or\\\n(ii)
$\g$ is an abelian extension of $\fa_3$.

\nit{Proof} We shall exhibit an orthogonal transformation of $(\g,\alr)$ which
reverses the orientation of $\g$ in each case. This automorphism will induce
the desired isometry of $(G,g)$.

\n(i) If $\fz\ne0$, choose an orthonormal basis $(\f_1,\f_2,\f_3,\f_4)$ of $\g$
with $\f_1\in\fz$. It is then immediate that the linear mapping
$\phi\colon\g\to\g$ defined by \[ \phi(\f_i)=\left\{\ba{ll}-\f_1\qq&i=1,\\
\f_i& i\ne1,\ea\right.\] is a Lie algebra automorphism.

\n(ii) If $\g =\R\op_\rho\fa_3$, choose an orthonormal basis $(\f_1,\ldots,\f_4)$
with $\f_2,\f_3,\f_4\in\fa_3$. This time, \[\psi(\f_i)=\left\{\ba{ll} \f_1\qq&
i=4,\\-\f_i\q& i\ge2\ea\right.\] is the required automorphism.\qed

Combined with Proposition~2.1, this yields

\nbf{Corollary} Let $G$ be an oriented Lie group with a left-invariant
Riemannian metric that satisfies $W^+=0$ and $W^-\ne0$. Then $\g$ is an abelian
extension of $\fa_2$ or $\fh_3$.

\nbf{Proposition} The Lie algebras $\g_\t,\,\g_{\t'}$ are isomorphic if
and only if $\t=\pm\t'$.

\nit{Proof} The derived algebra $\g_\t'=[\g_\t,\g_\t]$ in \rf{LAt} has a basis
$(\f_2,\f_3,\f_4)$. With respect to this basis, $\ad\,\f_1$ acts on $\g_\t'$ as \[
\left(\ba{ccc} 1&\t&0\\-\t&1&0\\0&0&2\ea\right).\] It follows that if $f$ is
any element of $\g_\t\setminus\g_\t'$ then
\[\frac{\det(\ad(f)|_{\g_\t'})}{\tr(\ad(f)|_{\g_\t'})}=\ft12(1+\t^2),\]
so $\g_\t,\g_{\t'}$ cannot be isomorphic unless $\t=\pm\t'$. We observed in \S1
that they \textit{are} isomorphic if $\t=-\t'$.\qed

\nbf{Proposition} Let $\t\in\R$ and $\s,\s'>0$. The Riemannian manifolds
$(G_\t,g_\s)$ and $(G_{\t'},g_{\s'})$ are isometric if and only if $\s=\s'$.

\nit{Proof} The fact that $(G_\t,g_\s)$ is not isometric to $(G_{\t'},g_{\s'})$
if $\s\ne\s'$ follows immediately from {\bf1.5}, the computations of which can
be carried out by starting from \rf{dei} below. It therefore suffices to show
that $(G_\t,g_\s)$ is isometric to $(G_0,g_\s)$ for fixed $\s>0$.

Define $\fu$ as the abelian extension \[ \fu=\al \f_0\ar\op_\rho \g_0,\] where
$\rho$ is given by \be\ba{ll} [\f_0,\f_2]=-\f_3,\qq&[\f_0,\f_1] = 0\\{} [\f_0,\f_3]=
\f_2,& [\f_0, \f_4]= 0.\ea\ee{LAu} Let $U$ denote the simply connected Lie group
with Lie algebra $\fu$.

Comparing \rf{LAu} with \rf{LAt}, we see that the mapping $\g_\t\to\fu$ defined by
\[\ba{rcl}\f_1&\mapsto& \t \f_0+\f_1,\\ \f_i &\mapsto& \f_i,\q 2\le
i\le4,\ea\] is a Lie algebra homomorphism. We denote the induced immersion
$G_\t\to U$ by $i_\t$.

The 1-parameter subgroup $\R$ of $U$ generated by $\f_0$ acts on the coordinates
$x,y$ by \be\ba{rcl} x &\mapsto& x\cos s+y\sin s,\y y &\mapsto&\!-x\sin s+y\cos
s.\ea\ee{cis} This gives rise to a submersion $\pi:U\to U/\R$ onto a
homogeneous manifold for which $\pi\cir i_\t$ is a diffeomorphism for each
$\t$. By construction, \[(\pi\cir i_\t)_*\f_i=(\pi\cir i_0)_*\f_i,\q 1\le i\le
4,\] and $\phi=(\pi\cir i_\t)^{-1}\cir(\pi\cir i_0)$ is the required isometry
$(G_0,g_\s)\to(G_\t,g_\s)$

More explicitly, we may identify $U$ with $\R\times G_0$ by means of
coordinates $s,t,x,y,u$ relative to which $\f^0=ds$ and \rf{1f} holds with
$\t=0$. Then \rf{cis} and the fact that $i_\t^*(ds)=\t dt$ imply that
\[\ba{l}\phi^*(dx^2+dy^2)=(dx+\t ydt)^2+(dy-\t xdt)^2,\y
\phi^*(du+\ft12(xdy-ydx))=du+\ft12(xdy-ydx)-\ft12\t(x^2\!+\!y^2)dt,\ea\]
whence $\phi^*g_\s=g_\s$.\qed

\setcounter{equation}0\setcounter{subsection}3
\subsection*{3. Proof of the theorem}

Let $G$ be a 4-dimensional Lie group with a chosen orientation, and let $g$ be
a left-invariant Riemannian metric. The latter corresponds to an inner product
$\alr$ on the Lie algebra $\g$ of $G$. The curvature tensor of $g$ and its
components $W^\pm$ are completely determined by the structure constants of the
Lie algebra $\g$ and its inner product, so the problem is purely algebraic.

Because we are only interested in metrics that are not conformally flat, the
Corollary tells us that $\g$ must be an abelian extension of $\fa_2$ or
$\fh_3$. This information alone is sufficient to proceed with the calculations,
though these are simplified by appealing to a classification of 4-dimensional
solvable Lie algebras that appears for example in \cite{Ver}. The algebras we
need to consider are listed in the following table, with the structure
constants encapsulated in the exterior derivatives of a basis
$(\f^1,\f^2,\f^3,\f^4)$ of $\g^*$. In each case, $d\f^1=0$. The real parameters
$\a,\b$ satisfy $0\le\a\le2$ and $\b\ge0$.

\bigbreak

\[\ba{|c|c|c|} \hline\hbox{type}\ph & (d\f^2,d\f^3,d\f^4)
&\hbox{extension of}\y\hline \g_2\op\g_2\ph & (0,\f\,^{13},\,\f^{24}) & \fa_2\y
\g_{4,1} & (0,\,\f^{13},\,\f^{14}+\f^{23}) & \fa_2,\ \fh_3\y 
\g_{4,2} & (0,\,\f^{13}-\f^{42},\,\f^{14}-\f^{23}) & \fa_2\y 
\g_{4,9}(\a) & ((1-\a)\f^{12},\,-\f^{13},\,-\a \f^{14}-\f^{23}) & \fh_3\y 
\g_{4,10} & (\f^{12},\,\f^{12}+\f^{13},\,\f^{23}+2\f^{14}) & \fh_3\y 
\g_{4,11}(\b) &\q(\b \f^{12}+\f^{13},\,-\f^{12}+\b \f^{13},\,-\f^{23}+2\b \f^{14})
\q&\fh_3\y\hline\ea\]

\vskip10pt\bigbreak

Let $\g$ be one of $\g_{4,1}$, $\g_{4,9}(\a)$, $\g_{4,10}$, $\g_{4,11}(\b)$.
In each case, the above basis satisfies \be\left.\ba{rcl} d\f^1\=0,\y
d\f^2&\in&\al \f^{12},\f^{13}\ar,\y d\f^3&\in&\al \f^{12},\f^{13}\ar,\y d\f^4&\in&\al
\f^{23},\f^{14}\ar.\ea\right.\ee{inc} Apply the Gram-Schimdt process to obtain an
orthonomal basis $(\e^1,\e^2,\e^3,\e^4)$ satisfying \be \e^j=\sum_{i=1}^j
a_i^j\f^i,\q a_j^j\ne0.\ee{GS} The relations \rf{inc} then become
\be\left.\ba{rcl} d\e^1\=0,\y d\e^2\= c^2_{12}\e^{12}+c^2_{13}\e^{13},\y d\e^3\=
c^3_{12}\e^{12}+c^3_{13}\e^{13},\y d\e^4\= c^4_{12}\e^{12}+c^4_{13}\e^{13}+
c^4_{14}\e^{14}+c^4_{23}\e^{23}, \ea\right.\ee{on} for certain real constants
$c_{ij}^k$.

In the new scheme of things, each line on the right-hand side must include all
terms from the previous line. Even if the isomorphism class of $\g$ is fixed,
the structure constants $c^k_{ij}$ are allowed to vary to reflect all possible
choices of inner product. However, we must first ensure that \rf{on} still
describes a Lie algebra. 

\nbf{Lemma} Further to \rf{on}, we may suppose that (i) $c^4_{23}=1$, (ii)
$c^4_{14}=c^2_{12}+c^3_{13}$, and (iii) $c^4_{12}=0$.

\nit{Proof} The Jacobi identity is equivalent to the set of equations
$\{d(d\e^i)=0:1\le i\le4\}$. In our case, only the equation $d(d\e^4)=0$ is in
doubt, and this yields \[ 0=c^4_{23}(c^2_{12}+c^3_{13}-c^4_{14}).\] For the Lie
algebras under consideration, $\f^{23}$ appears with a non-zero coefficient in
$d\f^4$. It follows from \rf{GS} that $c^4_{23}\ne0$ whence (ii). 

Given the similar nature of $d\e^2,d\e^3$, the form of \rf{on} is preserved by a
rotation of the form \be\ba{ll} \e^2\ \mapsto&\>(\cos\th)\e^2+(\sin\th)\e^3,\y
\e^3\ \mapsto&\!-\!(\sin\th)\e^2+(\cos\th)\e^3.\ea\ee{rot} Applying this
substitution, \[d\e^4=(c^4_{12}\cos\th-c^4_{13}\sin\th)\e^{12}+\cdots\] and we
can choose $\th$ so that (iii) is satisfied relative to the new orthonormal
basis.

Finally, performing the overall scaling $\e^i\mapsto c^4_{23}\e^i$ gives (i).
\qed

The Weyl tensor can now be computed in terms of the $c^k_{ij}$. The orientation
of $(\e^i)$ can be reversed by changing the sign of $\e^1$. This will alter the
signs of the $c^k_{1j}$, but will not interfere with the above assumptions.  It
therefore suffices to compute $W^+$, whose vanishing is not affected by the
rescaling that accompanied condition (i).

We now proceed to solve $W^+=0$, that represents $\dim S^2_0(\E^+\g^*)=5$
quadratic equations in the coefficients of \rf{on}. There remain 5 unknowns,
namely \[c^2_{12},\ c^2_{13},\ c^3_{12},\ c^3_{13},\ c^4_{13}.\]

\nbf{Proposition} Given the above assumptions, there are only two real
solutions of $W^+=0$, namely
\vspace{-8pt}\[\ba{ll}\hbox{(i)}& a^3_{13}=-1,\ a^4_{13}=0,\ 
a^3_{12}=-a^2_{13},\ a^2_{12}=-1;\\[6pt] \hbox{(ii)}& a^3_{13}=-\ft12,\
a^4_{13}=0,\ a^3_{12}=-a^2_{13},\ a^2_{12}=-\ft12.\ea\]

\nit{Proof} Involves an application of the Appendix.\qed
 
Setting $a^2_{13}=\t/\s$ with $\s=1$ or 2 gives \be \left.\!\ba{rcl}\ds
d\e^1\=0,\y\s\,d\e^2\=-\e^{12}-\t \e^{13},\y\s\,d\e^3\=\>\t \e^{12}-\e^{13},\y
\s\,d\e^4\=-2\e^{14}+\s \e^{23}.\ea\right.\ee{dei} Let $(\e_i)$ denote the dual
basis. Then $(\f_1,\f_2,\f_3,\f_4)=(\s \e_1,\e_2,\e_3,\e_4)$ satisfies \rf{dfi}
and \rf{ip} (the new $f$'s are unrelated to those in \rf{GS}). In terms of the
table, it is easy to verify that
\[\g_\t\cong\left\{\ba{ll}\g_{4,9}(2),\q&\hbox{if }\t=0,\y
\g_{4,11}(\frs1\t),&\hbox{if }\t>0.\ea\right.\]

To complete the proof of the theorem, similar arguments are required
for $\g_2\op\g_2$ and $\g_{4,2}$. In the first case, the analogue of
\rf{on} is \[\left.\ba{rcl} d\e^1\=0,\y d\e^2\=0,\y 
d\e^3\= c^3_{12}\e^{12}+c^3_{13}\e^{13},\y
d\e^4\= c^4_{12}\e^{12}+c^4_{13}\e^{13}+c^4_{14}\e^{14}+ c^4_{23}\e^{23}
+c^4_{24}\e^{24}.\ea\right.\] The Jacobi identity implies that 
\[ a^4_{23}(a^3_{13}-a^4_{14})+a^4_{13}a^4_{24}=0.\] An overall scaling
can be used to set $a^4_{24}=1$ and this leaves 5 unknowns. There are no
solutions to $W^+=0$.

The equations for $\g_{4,2}$ are initially more complicated. Gram-Schmidt
yields \[\left.\ba{rcl} d\e^1\=0,\y d\e^2\=0,\y d\e^3\=c^3_{12}\e^{12}+
c^3_{13}\e^{13}+c^3_{14}\e^{14}+c^3_{23}\e^{23}+c^3_{24}\e^{24},\y d\e^4\=
c^4_{12}\e^{12}+c^4_{13}\e^{13}+c^4_{14}\e^{14}+c^4_{23}\e^{23}+c^4_{24}\e^{24}.
\ea\right.\] A rotation in the plane generated by $\e^1,\e^2$, in analogy to
\rf{rot}, enables us to arrange that $c^3_{14}=0$ and $c^3_{24}\ne0$. After an
overall scaling of the basis we may assume that $c^3_{24}=1$.  The Jacobi
identity then implies that \[c^4_{13}=0,\q c^4_{13}=c^3_{13}.\] This leaves 6
unknowns, sufficiently few for the program to reveal that once again there are
no solutions with $W^+=0$.\mb

\setcounter{equation}0\setcounter{subsection}4
\subsection*{4. Further properties}

\n\textbf{Invariant complex structures.} To check whether the manifold
$(G_\t,g_\s)$ is hyperhermitian, we take a standard basis \[\ba{l}\w_1=
\e^{12}+\e^{34},\y\w_2=\e^{13}+\e^{42},\y\w_3=\e^{14}+\e^{23} \ea\] of
self-dual 2-forms. Let $I_1,I_2,I_3$ denote the corresponding left-invariant
almost complex structures, and $\th_1,\th_2,\th_3$ the Lee forms defined by
$d\w_i=\w_i\we\th_i$. An easy consequence of \rf{dei} is that \[\ba{l}\ds
\th_1=-\frac3\s \e^1-\t \e^4=\th_2,\\[8pt]\ds\th_3=-(1+\frac2\s)\e^1.\ea\] The
equality between $\th_1$ and $\th_2$ implies (via some standard theory 
\cite{S}) that $I_3$ is always integrable. On the other hand, $\th_1$ and
$\th_3$ are only equal if $\t=0$ and $\s=1$.

\nbf{Corollary} The Lie group $G_\t$ admits an invariant complex structure
$I_3$ for any $\t$, but an invariant hypercomplex structure $(I_1,I_2,I_3)$
only for $\t=0$. This hypercomplex structure is compatible with $g_1$.\mb

\bb\n\textbf{Properties of \hbox{\boldmath $g_1$}.} Suppose that $\t=0$ and
$\s=1$.  Using the coordinates $t,x,y,u$ as in \rf{1f} and the substitution
$v=\e^t$, take \[\e^1=dt,\q\e^2=\frac1v dx,\q \e^3=\frac1v
dy,\q\e^4=\frac1{v^2}(du+\ft12(xdy-ydx)).\] If we further set \[\ba{rcl}
\tw_1^\pm\=vdv\we dx\pm(dy \we du+\ft12ydy\we dx)\\[6pt]\tw_2^\pm\=vdv\we
dy\pm(du\we dx+\ft12xdx\we dy)\\[6pt]\tw_3^\pm\=dv\we du-\ft12(xdy\we dv
-ydx\we dv)\pm vdxdy,\ea\] then $\tw_i^+=v^3\w_i$. The following is immediate.

\nbf{Lemma} Five of the six 2-forms $\tw_i^\pm$ are closed, all except
$\tw_3^-$.\mb

The fact that $d\tw_i^+=0$ for all $i$ implies that the almost-complex
structures $I_1,I_2,I_3$ are all integrable, and $v^3g_1$ is a hyperk\"ahler
metric. This observation is contained in \cite{Bar2}; in fact \be v^3g_1=
v(dx^2+dy^2+dv^2)+\frac1v(du+\ft12(xdy-ydx))^2\ee{GH} is in the form of the
Gibbons-Hawking ansatz for a hyperk\"ahler metric with $S^1$ action. On the
other hand, if $\sigma=(\cos\th)\w_1^-+(\sin\th)\w^-_2$ then (for any $\th$)
$(v^3g_1,\sigma)$ is an almost-K\"ahler Einstein structure which is not
K\"ahler. There is a family of similar local counterexamples to the Goldberg
conjecture \cite{Arm}.

\bb\n\textbf{Links with other geometries.} Let $z=x+iv$, and consider the upper
half plane $H=\{z\in\C:v>0\}$. Then $v^3g_1$ coincides with the natural metric
on the cotangent bundle $T^*H$ induced from a `special K\"ahler' metric on
$H$. This aspect of the above example is explained in \cite{S}. Here, we merely
remark that special K\"ahler metrics are determined locally by a `holomorphic
prepotential' $\cF$ that here we take to be $\ft16z^3$ \cite{Fre}. This allows
one to introduce a dual or `conjugate' holomorphic coordinate \[w=\ft12z^2=
\ft12(x^2-v^2)-xv\,i,\] whose real part we denote by $\tx$. The (incomplete)
special K\"ahler metric $v(dx^2+dv^2)$ on $H$ has associated 2-form
\[\w=dx\we d\tx=\frac v{2i}dz\we d\overline{z}.\] Decreeing the 1-forms
$dx,d\tx$ to be parallel defines a flat symplectic connection on $H$ which can
be used to construct the hyperk\"ahler metric \rf{GH} in a functorial manner.

The Lie algebra $\g_0$ also features in \cite{Gib}, as one example in a study
of Einstein metrics constructed from left-invariant metrics on nilpotent
groups. Higher dimensional examples give rise to Einstein metrics on solvable
extensions, and related metrics with exceptional holonomy.

\bb\n\textbf{Quadratic components of \textit{W}.} In four dimensions,
the space $\cW$ of Weyl tensors has symmetric square \be S^2\cW\cong
S^2\cW^+\op(\cW^+\ot\cW^-)\op S^2\cW^-,\ee{7} and there exist
$SO(4)$-irreducible 9-dimensional representations $\cV^\pm$ for which
\[ S^2\cW^\pm\cong\cV^\pm\op \cW^\pm\op\R\] (see for example \cite{Sh4}). 
Thus, the vanishing of $\cW^+$ reduces the number of non-zero components of
$W\ot W$ in \rf{7} from 7 to 3. Moreover, if $W^+$ has repeated eigenvalues
as in {\bf1.5} then the $\cW^+$-component of $W^+\ot W^+$ is proportional to
$W^+$ (this also implies that certain components of $W^+\ot W^+\ot W^+$ vanish
in $S^3\cW^+$).

On a higher dimensional manifold $M$, it is known that the Weyl tensor $W$ is
irreducible. Appropriate analogues of anti-self-duality can in theory be
defined in terms of components of $S^2\cW$. Of special interest are likely to
be conditions imposed on $W$ by the existence of one or more orthogonal complex
structures on $M$.

\setcounter{equation}0\setcounter{subsection}5
\subsection*{5. Appendix}

The essential part of the program used to solve the equation $W^+=0$ is
reproduced below. It is set up to mimic the assumptions \rf{on} and {\bf3.4},
through the remaining cases of $\g_2\op\g_2$ and $\g_{4,2}$ are similar.

Once the solutions are found it is an easy matter to substitute back using for
example \texttt{assign(so[1])} and determine $W^-$ by setting \texttt{z:=-1}.

\begin{verbatim}
# Define the structure constants
   for i to 4 do defform(e[i]=1) od:
   for i to 5 do defform(a[i]=const) od:
   defform(d(e[1])=0):
   defform(d(e[2])=a[1]*e[1]&^e[2]+a[2]*e[1]&^e[3]):
   defform(d(e[3])=a[3]*e[1]&^e[2]+a[4]*e[1]&^e[3]):
   defform(d(e[4])=a[5]*e[1]&^e[3]+(a[1]+a[4])*e[1]&^e[4]+e[2]&^e[3]):

# Define the set of connection coefficients
   gset:={}:
   for k to 4 do g[k]:=array(antisymmetric,1..4,1..4) od:
      for k to 4 do for j to 4 do for i to j-1 do 
         defform(g[k][i,j]=const); gset:=gset union {g[k][i,j]}
      od od od:

# Define and solve the Cartan equations
   fset:={}: 
   for k to 4 do 
      de[k]:=sum(sum(g[m][k,n]*e[m]&^e[n],m=1..4),n=1..4) 
   od:
   c:=(k,i,j)->coeff(de[k]-d(e[k]),e[i]&^e[j]):    
   for k to 4 do for j to 4 do for i to j-1 do 
      fset:=fset union {c(k,i,j)=c(k,j,i)} 
   od od od:
   sol:=solve(fset,gset);
   assign(sol):

# Define the Riemann curvature tensor
   r1:=(i,j)->sum(g[p][i,j]*d(e[p]),p=1..4):
   r2:=(i,j)->sum(sum(sum(g[u][i,r]*g[v][r,j]*e[u]&^e[v],u=1..4),
                                                 v=1..4),r=1..4):
   r:=(i,j,x,y)->coeff(simpform(r1(i,j)-r2(i,j)),e[x]&^e[y]):
   for i to 4 do for j to 4 do
   R[i,j]:= array(antisymmetric,1..4,1..4):
      for y to 4 do 
         for x to y-1 do R[i,j][x,y]:=r(i,j,x,y)-r(i,j,y,x) od 
   od:
   od od:

# Define the Ricci, scalar and Weyl+ curvature
   ric:=array(symmetric,1..4,1..4):
   for j to 4 do for i to j do 
      ric[i,j]:= sum(R[i,p][j,p],p=1..4) 
   od od:
   sc:=sum(ric[q,q],q=1..4):
   W:=array(symmetric,1..3,1..3): z:=1:
      W[1,1]:= R[1,2][1,2]+R[3,4][3,4]+2*z*R[1,2][3,4]-sc/6:
      W[2,2]:= R[1,3][1,3]+R[4,2][4,2]+2*z*R[1,3][4,2]-sc/6:
      W[3,3]:= R[1,4][1,4]+R[2,3][2,3]+2*z*R[1,4][2,3]-sc/6:
      W[1,2]:= R[1,2][1,3]+z*R[1,2][4,2]+z*R[3,4][1,3]+R[3,4][4,2]:
      W[1,3]:= R[1,2][1,4]+z*R[1,2][2,3]+z*R[3,4][1,4]+R[3,4][2,3]:
      W[2,3]:= R[1,3][1,4]+z*R[1,3][2,3]+z*R[4,2][1,4]+R[4,2][2,3]:

# Solve Weyl+=0
      hset:={W[1,1]=0,W[1,2]=0,W[1,3]=0,W[2,2]=0,W[2,3]=0}:
      so:=solve(hset):
\end{verbatim}

\renewcommand{\thebibliography}{\list{{\bf\arabic{enumi}.\hfil}}
{\settowidth\labelwidth{18pt}\leftmargin\labelwidth\advance
\leftmargin\labelsep\usecounter{enumi}}\def\newblock{\hskip.05em} \sloppy
\sfcode`\.=1000\relax} \newcommand{\bi}{\vspace{-4pt}\bibitem}

\small\parskip10pt

\n Rue Emile Bouilliot 25, 1050 Bruxelles, Belgium

\n Dipartimento di Matematica, Politecnico di Torino, 10129 Torino, Italy

\enddocument